\documentclass[a4paper]{article}
\usepackage{amsmath}
\usepackage{amssymb}
\usepackage{amsthm}

\newtheorem{thm}{Theorem}[section]
\newtheorem{cor}[thm]{Corollary}
\newtheorem{lem}[thm]{Lemma}
\newtheorem{prop}[thm]{Proposition}
\theoremstyle{definition}
\newtheorem{defn}[thm]{Definition}
\theoremstyle{remark}
\newtheorem{rem}[thm]{Remark}
\numberwithin{equation}{section}
\newcommand{\norm}[1]{\left\Vert#1\right\Vert}
\newcommand{\abs}[1]{\left\vert#1\right\vert}
\newcommand{\set}[1]{\left\{#1\right\}}
\newcommand{\Real}{\mathbb R}
\newcommand{\eps}{\varepsilon}

\newcommand{\para}[1]{\left(#1\right)}
\newcommand{\inp}[1]{\left\langle#1\right\rangle}

\title{ Large Deviation Principle for Semilinear Stochastic Evolution Equations with Monotone Nonlinearity and Multiplicative Noise}%
\author{{Hassan Dadashi-Arani  and   Bijan Z. Zangeneh}}
\begin{document}

\date{}%

\maketitle
\begin{center}
\emph{Department of Mathematical Science, Sharif University of Technology, Tehran, Iran}
\end{center}
\let\thefootnote\relax\footnotetext{
\hspace{-.5cm}\textbf{Subject Class:}{ 60F10, 60H20}
This work is partially supported by the research council of Sharif University of Technology.}

\begin{abstract}
We demonstrate the large deviation property for the mild solutions of stochastic evolution equations with monotone nonlinearity and multiplicative noise. This is achieved using the recently developed weak convergence method, in studying the large deviation principle.
An It\^{o}-type inequality is a main tool in the proofs. We also give two examples to illustrate the applications of the theorems.
\end{abstract}
\section{\textbf{Introduction}}
\begin{table}[b]
\begin{tabular}{l}
 \hline

\end{tabular}
\end{table}

The limit theorems of probability theory, such as the strong law
of large numbers and the central limit theorem, basically say
that averages taken over large samples converge to expected
values. These results say little or nothing about the rate of
convergence which is important for many applications of
probability theory, e.g., statistical mechanics. One way to
address this is the theory of large deviations. Some basic ideas
of the theory can be tracked back to Laplace and Cram\'{e}r,
although a clear unified formal definition was introduced
by Varadhan in \cite{kn:Va}. Varadhan formulated the \emph{large
deviation principle} (LDP) for a collection of random variables
$\{X^{\eps}: \eps>0\}$ on a Polish space $\mathcal{E}$, a
complete separable metric space.
\begin{defn}
Let $(\Omega, \mathcal{F}, \textbf{P})$ be a probability space. We say that the
$\mathcal{E}$-valued family of random variables $\{X^{\eps}: \eps>0\}$ obeys the LDP with the good \emph{rate function} (or \emph{action functional}) $I$, where $\emph{I}:\mathcal{E}\rightarrow [0,\infty]$ is a lower semicontinuous function, if\\

$A1)$ the level set $K_N:=\{x\in \mathcal{E} : \emph{I}(x)\leq N\}$ is compact for every $N< \infty$;\\

$A2)$ for any open set $G$ in $\mathcal{E}$,
\begin{align}
\limsup\limits_{\eps\rightarrow0} \eps^2 \log \textbf{P}(X^{\eps}\in G)\leq -\inf\limits_{x\in G}\emph{I}(x);\label{33}
\end{align}

$A3)$ for any closed set $F$ in $\mathcal{E}$,
\begin{align}
\liminf\limits_{\eps\rightarrow0} \eps^2 \log \textbf{P}(X^{\eps}\in F)\geq -\inf\limits_{x\in F}\emph{I}(x).\label{9}
\end{align}
\end{defn}
In the literature of large deviations for stochastic differential equations, $\mathcal{E}$ is the space of solutions and $X^{\eps}$ is the solution, when the intensity of noise is multiplied by $\eps$. To establish the LDP for solutions of stochastic differential equations, the following lemma is a crucial tool.

\begin{lem} (\textbf{Varadhan's Contraction Principle}) Let $\{X^{\eps}, \eps>0\}$ be a family of random variables on the Polish space $\mathcal{E}$ which satisfies the LDP with the rate function $I:\mathcal{E}\rightarrow\Real$, and $f:\mathcal{E}\rightarrow \mathcal{E}^\prime$ be a continuous transformation. Then the family $\{f(X^{\eps}), \eps>0\}$ satisfies the LDP on the Polish space $\mathcal{E}^\prime$ with the rate function
$$I^\prime(y)=\inf\limits_{\{x \in \mathcal{E} : \; f(x)=y\}} I(x),\quad y\in \mathcal{E}^\prime.$$
Here we obey the usual convention $\inf \emptyset = +\infty$.
\end{lem}
Using this lemma and Schilder's theorem \cite{kn:Sc}, we can immediately establish the LDP for SDE's with additive noise. In their pioneering book \cite{kn:FW}, Freidlin and Wentzell used a time discretization argument to freeze the diffusion term in the case of SDE's with multiplicative noise. In each step of time discretization process, the solution is a continuous function of the noise and therefore satisfies the LDP by the Varadhan's contraction principle. Finally, one must prove that the equation with the frozen drift is a good exponential estimate of the original solution. To apply this method in infinite dimensions, after freezing the diffusion, one must project it to a finite dimensional system. Then, after establishing the LDP in finite dimension, one must further argue that the LDP remains valid as one approaches the infinite dimensional system. One can find the fundamental ideas of this method for infinite dimensional spaces in \cite{kn:AZ, kn:PR, kn:Le}. We refer to Freidlin \cite{kn:Fr} for a first work in studying LDP for stochastic reaction-diffusion equations. Thereafter, many papers in different approaches to SPDE have been written in any of which various technical difficulties have been overcome. In his thesis \cite{kn:Pe} and in \cite{kn:Pes}, Peszat used these fundamental ideas to establish LDP for semilinear stochastic evolution equations with Lipschitz nonlinearity and multiplicative noise. As well, Sowers in \cite{kn:So} and \cite{kn:Sow} studied the LDP problem for stochastic reaction-diffusion equations with diffusion coefficients bounded on the unit circle. Cerrai and R\"{o}ckner in \cite{kn:CR} investigated the large deviation property for stochastic reaction–-diffusion systems with unbounded diffusion terms and locally Lipschitz reaction terms which satisfy a polynomial growth.

In this paper, we follow a totally different (weak convergence) approach to demonstrate the LDP. In this approach, which has been initiated recently in \cite{kn:BD, kn:BDM}, we consider an equivalent formulation of the LDP, called the Laplace principle. To this end we must consider the right hand sides of (\ref{33}) and ($\ref{9}$) for bounded continuous functions instead of indicators of closed and open sets. See \cite[Section 1.2]{kn:DE} for a proof of equivalence between Laplace principle and large deviation principle.

\begin{defn} (\textbf{Laplace principle}) The family of random variables $\{X^\eps\}$ on the Polish space $\mathcal{E}$ is said to satisfy the \emph{Laplace principle} with the \emph{rate function I} if for all bounded continuous functions $h:\mathcal{E}\rightarrow\Real$,
$$\lim_{\eps\rightarrow0}\eps^2 \log \mathbb{E}\;{\exp\left\{-\frac{1}{\eps^2}h(X^\eps)\right\}}=-\inf\limits_{f\in \mathcal{E}}\{h(f)+\emph{I}(f)\}.$$
\end{defn}
The weak convergence approach is ideally suited for evaluation of the integrals appearing in the Laplace principle. Briefly, the basic setup in this approach is to associate with the given LDP problem a family of minimal cost functions which gives a variational representation for these integrals. Then the asymptotics of the minimal cost functions will be determined by the weak convergence methods. Dupuis and Ellis in \cite{kn:DE} present a comprehensive reference of weak convergence methods. For some illustrations of this method in large deviations, we refer to \cite{kn:Li, kn:RZZ, kn:SS}. Liu  \cite{kn:Li} deals with a problem similar to the one considered in the present paper in the variational approach to SPDE.
Knowing that some physical (such as hyperbolic) systems can not be treated by this approach, we show the LDP or hyperbolic systems as an application of our main theorem. To state the main theorem of weak convergence method, we need a short preliminary.

Let  $H$ be a real separable Hilbert space with the norm and inner product denoted by $\inp{\cdot,\cdot}$ and $\norm{\cdot} $, respectively. Assume that $(\Omega,\mathcal{F},\mathcal{F}_t,\textbf{P})$ is a complete stochastic basis with a right-continuous filteration $\{\mathcal{F}_t, t\geq 0\}$. Let $\{e_k\}$ be an orthonormal basis for $H$ and $\{W_k\}$ be a sequence of independent, real-valued, $\mathcal{F}_t$-Wiener processes. By a \emph{cylindrical Wiener process} on $H$ we mean the series
$$W(t)=\sum\limits_{k\in \mathbb{N}} W_k(t) e_k, \quad t\geq 0.$$
This series does not converge in $H$ but it converges in an arbitrary Hilbert space $U$ containing $H$ with a Hilbert-Schmidt embedding. We denote the norm on $U$ by $\norm{\cdot}_U$. A function $f:[0,T] \times \Omega\rightarrow \Real$ is called  \emph{elementary} if there exist $a,b\in [0,T], \; a\leq b$, and a bounded $\mathcal{F}_a$-measurable random variable $X$ such that $f(s,\omega)=X(\omega) 1_{(a,b]}(s)$. Any finite sum of elementary functions is referred to as a simple function. We denote by $\bar{\mathcal{S}}$ the class of all simple functions. The predictable $\sigma$-field $\mathcal{P}$ on $\Omega\times[0,T]$ is the $\sigma$-field generated by $\bar{\mathcal{S}}$. For a given Hilbert space $U$, a  $\mathcal{P}$-measurable function $u:\Omega\times[0,T]\rightarrow U$ is called \emph{$U$-valued predictable process}. Let $\mathcal{P}^2(U)$ be the family of all $U$-valued predictable processes $u$ for which $\int^T_0 \norm{u(s)}^2 ds < \infty\;$ a.s.

The following proposition is the first attempt to apply the  weak convergence theory in large deviations. It offers a different criterion in studying Laplace principle by making a connection between exponential functional and variational representation.
\begin{prop}\label{p1}
Let $h$ be a bounded, Borel measurable function mapping $\mathcal{C}([0,T];U)$ into $\Real$. Then for any cylindrical Brownian motion $W,$
\begin{align}
-\log \mathbb{E}\;\{\exp(-h(W))\}=\inf \limits_{u\in\mathcal{P}^2(U)}\mathbb{E}\;\para{\frac{1}{2}\displaystyle\int^T_0\norm{u(s)}^2_U ds+h\para{W+\displaystyle\int^._0u(s)ds}}\nonumber.
\end{align}
\end{prop}

Finally, we define for every $N>0$
$$\mathcal{S}_N(U):=\{u \in L^2([0,T];U) : \displaystyle\int^T_0\norm{u(s)}^2_U ds\leq N\},$$
$$\mathcal{P}_N(U):=\{u \in \mathcal{P}^2(U) : u(\omega)\in \mathcal{S}_N(U), \mathbb{P}-a.s.\}.$$\\
\textbf{Hypothesis 1.} There exist measurable maps $ \mathcal{G}^0$ and $\mathcal{G}^\eps:\mathcal{C}([0,T];U)\rightarrow \mathcal{C}([0,T];H)$, $\eps > 0$, such that the following two conditions hold:
\begin{description}
\item {(i)} for every $N<\infty$, the set
$$K_N:=\{\mathcal{G}^0(\displaystyle\int^._0u(s)ds):u\in \mathcal{S}_N(U)\}$$
is a compact subset of $\mathcal{C}\para{[0,T];H}$;
\item {(ii)} consider $N<\infty$ and the family $\{u^{\eps}, \eps>0\} \subset\mathcal{P}_N(U)$, such that $u^{\eps}$ converges to $u^0$ in the weak topology of $\mathcal{S}^N(U)$ and in distribution (as $\mathcal{S}_N(U)$-valued random variables), then
$$\mathcal{G}^{\eps}\para{\eps W+\displaystyle\int^._0u^{\eps}(s)ds}\rightarrow \mathcal{G}^0(\displaystyle\int^._0u(s)ds),$$
in the strong topology of $\mathcal{C}\para{[0,T];H}$ and in distribution (as $\mathcal{C}\para{[0,T];H}$-valued random variables), as $\eps \rightarrow 0$.
\end{description}
In this paper, the space  $\mathcal{C}([0,T];H)$ is assumed to be equipped with the (strong) topology of $\sup$-norm, unless stated explicitly.

The following theorem (Theorem 5 in \cite{kn:BDM}) is the main application of weak convergence theory in LDP.
\begin{thm}\label{h5}
Let $X^{\eps}=\mathcal{G}^{\eps}(\eps W)$, for every $\eps>0$, and assume Hypothesis 1 holds. For each $f\in \mathcal{C}([0,T];H)$, let
\begin{align}
I(f):=\inf \limits_{\{u\in L^2([0,T];U):f=\mathcal{G}^0(\int^._0u(s)ds)\}} \frac{1}{2}\int^T_0\norm{u(s)}^2_U ds,\label{22}
\end{align}
then the family $\{X^{\eps}:\;\eps>0\}$ satisfies the LDP with the good rate function $I$.
\end{thm}
The aim of this paper is to establish the LDP for the mild solution of semilinear stochastic evolution equation with monotone nonlinearity and multiplicative noise, using the above Theorem. In the next section, we express precise assumptions on the equation. When we apply the weak convergence method, the first step is to verify the well-posedness of the original deterministic equation controlled by a function $u \in L^2([0,T];U)$. This is called the skeleton problem and will be considered  in Section 3. In Section 4, the action functional will be investigated. Section 5 is devoted to verify Hypothesis 1(ii) which is the essential step when one applies the weak convergence method to demonstrate the LDP. Finally in Section 6, applications of the theorems in two examples will be discussed.

\section{\textbf{Assumptions and preliminaries}}

There are three different approaches dealing with the different problems in SPDE; namely, the variational approach, the martingale measure approach, and the semigroup approach. Krylov and Rozovskii\cite{kn:KR}, Walsh\cite{kn:Wa} and Da Prato and Zabczyk\cite{kn:DPZ} are the comprehensive references of these approaches, respectively. In this paper, we follow the semigroup approach for  semilinear stochastic evolution equations with monotone nonlinearities. One of the advantages of this approach comparing to the other approaches in SPDE is that it gives a unified treatment of a wide class of parabolic, hyperbolic and delay differential equations. The stressed assumption on the equation considered in this paper is the monotone nonlinearity. In some literature, especially in the numerical analysis literature, the monotonicity condition is also named the one-sided lipschitz condition. In \cite{kn:Z1, kn:Z2, kn:Za}, Zangeneh explored this notion for a general semilinear stochastic evolution equation extensively. He proved the existence, uniqueness, and some qualitative properties of the solution of the integral equation
\begin{align}
X_t=U(t,0)X_0+\int^t_0 U(t,s)\left(X_s+f(s,X_s)\right)ds+\int^t_0 U(t,s)g(s,X_s)dW_s+V(t),\label{16}
\end{align}
where $U(t,s),\; t\geq s\geq 0,$ is a strong semigroup with generator $A(t), \; t\geq 0$, $f$ is semimonotone, $g$ is Hilbert-Schmidt valued and Lipschitz, and $V$ is a continuous adapted process. The LDP for the solution of  $(\ref{16})$ with additive noise, when $g$ is constant, has been considered in \cite{kn:Zang}. In the case of $V=0$, the existence, uniqueness and stability of the continuous solutions of  $(\ref{16})$ with and without delay have been investigated in \cite{kn:JZ, kn:Ja2} and \cite{kn:Ja, kn:Ja1}, respectively. The existence, uniqueness and measurability of the solution of $(\ref{16})$ with a stopping time parameter have been studied in \cite{kn:HZ}. Another illustration of monotone nonlinearity in semilinear stochastic evolution equation appeared in \cite{kn:ZZ} in which the random motion of an elastic string is modeled when the elasticity tends to infinity.

To simplify notations, we assume that the linear operators $A(t), \; t\geq 0$, and the functions $f$ and $g$ are independent of $t$. Actually, in the absence of these restrictions, we only need a few trivial modification in our proofs. It is proved in \cite{kn:Z} that the solution $X$ is a continuous function of the adapted process $V$. Therefore, by Varadhan's contraction principle, this term does not give any serious difficulty in demonstrating the large deviation property. In this way, we consider the small noise LDP, as $\eps \rightarrow 0$, for the mild solution of the stochastic semilinear evolution equation
\begin{align}
dX^{\eps}_t=(AX^{\eps}_t+f(X^{\eps}_t))dt+\eps g(X^{\eps}_t)dW_t, \; \;  X(0)=X_0, \label{13}
\end{align}
which is, by definition, the strong solution of the  integral equation
\begin{align}
X^\eps_t=S(t)X_0+\displaystyle\int^t_0S(t-s)f(X^\eps_s)ds+\eps\displaystyle\int^t_0S(t-s)g(X^\eps_s)dW_s.\label{20}
\end{align}
The linear operator $A$, and the functions $f:H\rightarrow H$ and $g:H\rightarrow L_2(U,H)$ satisfy the following hypotheses, and $W$ is a cylindrical Brownian motion. We should remind that $L_2(U,H)$ and $L(U,H)$ are the spaces of Hilbert-Schmidt and bounded linear operators from $U$ into $H$ with the norms $\norm{\cdot}_2$ and $\norm{\cdot}_L$, respectively.\\\\
\textbf{Hypothesis 2}.
\begin{description}
  \item (i) $A:D(A)\rightarrow H$ is a closed linear operator with dense domain $D(A)\subseteq H$ which generates a $C_0$-semigroup $S(t)$ on $H$ and there exist $L, \lambda >0$ such that $\norm{S(t)}\leq L e^{\lambda t}$, for every $t\geq 0$;\\

  \item (ii) $f$ is a demicontinuous function on $H$, i.e., whenever $\{x_n\}$ converges strongly to $x$ in $H$, $\{f(x_n)\}$ converges weakly to $f(x)$ in $H$;\\

  \item (iii) $f$ satisfies a polynomial growth condition, i.e.
$$\exists m\in \mathbb{N}\; s.t. \quad \norm{f(x)} \leq C(1+\norm{x}^m ),\qquad \forall x\in H;$$

  \item (iv) $-f$ is semimonotone with parameter $M\geq 0$, i.e.
$$\inp{f(x)-f(y),x-y}\leq M\norm{x-y}^2 ,\qquad \forall x,y\in H;$$

  \item (v) $g$ is Lipschitz continuous with Lipschitz constant $D$, i.e.
$$\norm{g(x)-g(y)}_2\leq D \norm{x-y} ,\qquad \forall x,y\in H.$$
\end{description}

\begin{rem}
\begin{description}
\item (i) Setting $y=0$ in (iv), we obtain for every $x\in H$
$$\inp{f(x),x}\leq C(1+\norm{x}^2)$$

\item (ii) The immediate consequence of (v) is that $g$ satisfies a linear growth condition, i.e. $\norm{g(x)}_2\leq C(1+\norm{x})$. Moreover, since $\norm{\cdot}_L\leq \norm{\cdot}_2$, we have for all $x,y \in H$
\begin{align}
\begin{array}{c}
  \norm{g(x)-g(y)}_L\leq D \norm{x-y} \\ \\
  \norm{g(x)}_L\leq C(1+\norm{x}).\nonumber
\end{array}
\end{align}
\end{description}
\end{rem}
Now, we review some results of \cite{kn:Za, kn:Zan} which are regularly used  in the sequel.
\begin{prop}\label{p2}
Let $p\geq 2$. If  $ \mathbb{E}\norm{X_0}^p<\infty$ and Hypothesis 2 is satisfied, then $(\ref{13})$ has a unique continuous adapted mild solution $X^{\eps}$ such that
$$\sup\limits_{0<\eps\leq 1} \mathbb{E}\; (\sup\limits_{0\leq s\leq t} \norm{X^{\eps}_s}^p)< \infty, \quad \forall t\in [0,T].$$
\end{prop}
Actually this proposition has been proved in \cite{kn:Za}(Theorem 4) for solution to equation (\ref{16}) under some conditions similar to Hypothesis 2, and the proof has been stated only for $\eps=1$. But it is clear from the argument that it gives a uniform a priori bound for all $0< \eps\leq 1$.
\begin{prop}(\textbf{Energy-type inequality})
Let $a:[0,T]\rightarrow H$ be an integrable function. Suppose $A$ and $S$ satisfy Hypothesis 2(i). If
\begin{align}
X_t=S(t)X_0+\displaystyle\int^t_0S(t-s)a(s)ds,\nonumber
\end{align}
then
\begin{align}
\norm{X_t}^2\leq e^{2\lambda t}\norm{X_0}^2+2\displaystyle\int^t_0 e^{2\lambda(t-s)}\inp{X(s),a(s)}ds,\quad\forall t\in[0,T].\nonumber
\end{align}
\end{prop}
In Section 5, we use the stochastic version of this inequality obtained in \cite{kn:Za} which is given as
\begin{prop}(\textbf{It\^{o}-type inequality})
Let $\{Z_t: \; t\in[0,T]\}$ be an $H$-valued, cadlag, locally square integrable semimartingale. Suppose $A$ and $S$ satisfy Hypothesis 2(i). If
\begin{align}
X_t=S(t)X_0+\displaystyle\int^t_0S(t-s)dZ_s,\nonumber
\end{align}
then for all $t\in[0,T]$,
\begin{align}
\norm{X_t}^2\leq e^{2\lambda t}\norm{X_0}^2+2\displaystyle\int^t_0 e^{2\lambda(t-s)}\inp{X(s),dZ_s}+e^{2\lambda t}\left[\displaystyle\int^._0 e^{-\lambda s}dZ_s\right]_t,~~~w.p. 1,\nonumber
\end{align}
where $[ \quad ]_t$ stands for the quadratic variation process.
\end{prop}

The following lemma is a direct consequence of Burkholder's inequality for real martingales and Young inequality 
$$ab\leq K\frac{a^p}{p}+K^{-\frac{q}{p}}\frac{b^q}{q},$$ 
which is true for all positive real numbers $a, b, K$ and  $p,q$ satisfying $\frac{1}{p}+\frac{1}{q}=1$.
\begin{lem}\label{lem1} Let $X_t, \; t\in [0,\infty)$, be an $H$-valued continuous process. if $M_t$ is an $H$-valued continuous martingale, then for any constant $K>0$, we have
$$\mathbb{E}\left\{\sup\limits_{0\leq \theta\leq t}\abs{\displaystyle\int^{\theta}_0\inp{X_s,dM_s}}\right\}\leq\frac{3}{2K}\mathbb{E}(X^*_t)^2+\frac{3K}{2}\mathbb{E}([M]_t),$$
in which $X^{*}_t=\sup\limits_{0\leq s\leq t}\norm{X_s}$.
\end{lem}

To estimate the stochastic convolutions in Section 5, we need the following proposition.
\begin{prop}(\textbf{Burkholder-type inequality})
Suppose $A$ and $S$ satisfy Hypothesis 2(i), and $M$ is an $H$-valued continuous martingale. If $p\geq1$, then there exist suitable constant $C$ such that
\begin{align}
\mathbb{E}\left\{\sup\limits_{0\leq t\leq T}\norm{\displaystyle\int^t_0S(t-s)dM_s}^{2p}\right\}\leq Ce^{\lambda T}\mathbb{E}\{[M]^p_T\}, \quad \forall T>0.\label{18}
\end{align}\
\end{prop}
An extension of this proposition for stopping times has been proved in \cite{kn:HZ1}.

The following lemma, taken from \cite[Lemma 3.2]{kn:BD}, is extremely useful in our arguments.
\begin{lem}\label{lem2}
Let $\{u^n\}$ be a sequence in $\mathcal{P}_N(H)$. Suppose that $u^n$ converges in distribution to $u$ weakly in the space $L^2([0,T];U)$, then $\int^._0 u^n(t)dt$ converges in distribution to $\int^._0u(t)dt$ in the space $\mathcal{C}([0,T];H)$.
\end{lem}
In the case where $u^n$ is deterministic, this lemma is still meaningful and will be used in Section 4.

For the sake of completeness, we mention the well known  Gronwall inequality which is regularly used in our arguments.\\\\
\textbf{Gronwall's Lemma}. Let $f, \alpha, \beta :[0,T]\rightarrow \Real^+$ be a Lebesgue measurable functions such that $\beta$ is locally integrable and $\displaystyle\int^T_0\beta (s)f(s)ds<\infty$. Suppose
$$f(t)\leq \alpha(t)+\displaystyle\int^t_0\beta(s)f(s)ds, \; \; t\in [0,T],$$
then
$$f(t)\leq \alpha (t)+\displaystyle\int^t_0\alpha(s)\beta(s)\exp\left\{\int^t_s\beta(r)dr\right\}ds, \; \; t\in [0,T].$$
If  $\alpha$ is a constant function, then
$$f(t)\leq \alpha\exp\para{\displaystyle\int^t_0\beta(s)ds}.$$

\begin{rem} Since we investigate a limit problem on the finite interval $[0,T]$ and $\sup\limits_{0\leq t\leq T}\norm{S(t)}_L\leq e^{\lambda T}$, without loss of generality, we may assume that $S(t)$ is a semigroup of contractions.
\end{rem}

\section{\textbf{The skeleton problem}}
In this section, we study the following deterministic equation:
\begin{align}
\frac{d}{dt} z(t) = Az(t)+f(z(t))+g(z(t))u(t), \quad    z(0)=X_0,\label{14}
\end{align}
where $u \in L^2([0,T];U)$. We prove that (\ref{14}) has a mild solution and we obtain an a priori estimate for the solution in the space $\mathcal{C}\para{[0,T];H}$.
\begin{thm}\label{h1}
For each $u\in L^2([0,T];U)$ the integral equation
\begin{align}
z_t=S(t)X_0+\displaystyle\int^t_0S(t-s)\para{f(z_s)+g(z_s)u_s}ds,\label{15}
\end{align}
has a strong solution.
\begin{proof} First, we assume that $u \in L^\infty\para{\left[0 , T\right];U}$.
Browder \cite{kn:BR} and Kato \cite{kn:KA} showed the existence of mild solution for semilinear equations with semimonotone nonlinearities. So, it is sufficient to verify the semimonotonicity of
\begin{equation}
\begin{array}{l}
F:\Real \times H \mapsto H\\
F(t,x)=f(x)+g(x)u(t).
\end{array}\nonumber
\end{equation}
According to Hypothesis 2 (iv, v), we have for every $t \in \Real$ and $x, y\in H$
\begin{align}
\inp{F(t,x)-F(t,y),x-y}&=\inp{f(x)-f(y),x-y}+\inp{g(x)u(t)-g(y)u(t),x-y}\nonumber\\
&\leq\para{M+D\norm{u}_{\infty}}\norm{x-y}^2.\nonumber
\end{align}
So, $F$ is semimonotone with parameter $M^{\prime}=M+D\norm{u}_{\infty}$.
Then, since $L^{\infty}\para{[0,T];U}$ is dense in the space $ L^2\para{[0,T];U}$, an arbitrary $u\in L^2\para{[0,T];U}$ can be approximated by a sequence $u^n \in L^{\infty}\para{\left[0 , T\right]; U}, n\in \mathbb{N}$ in the strong topology of $L^2\para{[0,T];U}$. Let $z^n$ be the solution of  (\ref{15}) with the control function $u^n$. We show  that $\{z^n\}$ is a Cauchy sequence in $\mathcal{C}\para{\left[0 , T\right]; H}$. By the energy-type inequality, Hypothesis 2(iv), and the Cauchy-Schwarz inequality, we have
\begin{align}
&\norm{z^n_{t}-z^m_{t}}^2 \nonumber\\
&\leq2\displaystyle\int^t_0\inp{z^n_s-z^m_s,f(z^n_s)-f(z^m_s)+g(z^n_s)u^n_s-g(z^m_s)u^m_s}ds\nonumber\\
&\leq2M\int^t_0\norm{z^n_s-z^m_s}^2ds+2\int^t_0\inp{z^n_s-z^m_s,g(z^n_s)(u^n_s-u^m_s)}ds\nonumber\\
&\quad+2\int^t_0\inp{z^n_s-z^m_s,(g(z^n_s)-g(z^m_s))u^m_s}ds\nonumber\\
&\leq2M\int^t_0\norm{z^n_s-z^m_s}^2ds+2\int^t_0\norm{z^n_s-z^m_s}\norm{g(z^n_s)}_L\norm{u^n_s-u^m_s}_Uds\nonumber\\
&\quad+2\int^t_0\norm{z^n_s-z^m_s}\norm{g(z^n_s)-g(z^m_s)}_L\norm{u^m_s}_Uds.\label{39}
\end{align}
Applying the Young inequality to the second term and Remark 2.1(ii) to the third term, we obtain
\begin{align}
&\norm{z^n_{t}-z^m_{t}}^2 \nonumber\\
&\leq\int^t_0\left(2M+\norm{g(z^n_s)}^2_L+2D
\norm{u^m_s}_U\right)\norm{z^n_s-z^m_s}^2ds+\int^t_0\norm{u^n_s-u^m_s}^2_Uds.\label{3}
\end{align}
Now, we define $\alpha^{m,n}(t)=\displaystyle\int^t_0\norm{u^n_s-u^m_s}^2_Uds$  and  $\beta^{m,n}(t)=2M+\norm{g(z^n_t)}^2_L+2D\norm{u^m_t}_U$. Applying Gronwall's lemma, we get
\begin{align}
\norm{z^n_t-z^m_t}^2 \leq\alpha^{m,n}(t)+\displaystyle\int^t_0\alpha^{m,n}(s)\beta^{m,n}(s)
\exp\para{\displaystyle\int^t_s\beta^{m,n}(r)dr}ds.\label{4}
\end{align}
Noting  that $\{u^n\}$ is a Cauchy sequence in the strong topology of $L^2([0,T];U)$, we easily obtain
\begin{align}
\lim_{m,n\rightarrow\infty} \sup\limits_{0\leq t\leq T}\alpha^{m,n}(t)\rightarrow0.\label{31}
\end{align}
To estimate the first term in (\ref{3}) we need an a priori estimate for $z^n$. By an argument similar to that employed to derive (\ref{39}) and (\ref{3}), and using Remark 2.1(i), we have
\begin{align}
\norm{z^n_{t}}^2 &\leq \norm{X_0}^2+2\int^t_0\inp{z^n_s,f(z^n_s)+g(z^n_s)u^n_s}ds\nonumber\\
&\leq\norm{X_0}^2+2\int^t_0\para{C\norm{z^n_s}^2+\norm{z^n_s} \norm{g(z^n_s)}_L\norm{u^n_s}_U}ds\nonumber\\
&\leq \norm{X_0}^2+\int^t_0\para{\norm{z^n_s}^2 +\norm{f(z^n_s)}^2 +\norm{z^n_s}^2 \norm{u^n_s}^2_U+\norm{g(z^n_s)}^2_L}ds\nonumber\\
&\leq\displaystyle\int^t_0\left(\norm{u^n_s}^2_U+C\right)\norm{z^n_s}^2 ds+C. \label{5}
\end{align}
Note that, in this paper, $C$ will denote a positive constant whose exact value is not important and may change from line to line. Using Gronwall's lemma and a suitable constant $C$, we obtain the following estimate:
\begin{equation}
\norm{z^{n}_t}^2 \leq C\exp\para{\displaystyle\int^t_0\para{C+\norm{u^n_s}^2_U}ds},\nonumber
\end{equation}
which implies
\begin{align}
\sup\limits_{n\in \mathbb{N}}\sup\limits_{0\leq t\leq T}\norm{z^n_t} <\infty.\label{10}
\end{align}
Hence, we have
\begin{align}
&\displaystyle\int^t_0\beta^{m,n}(s)ds\nonumber\\
&=\displaystyle\int^t_0\para{2M+\norm{g(z^n_s)}^2_L+2D\norm{u^m_s}_U}ds\leq \int^t_0\para{C+C\norm{z^n_s}^2 +D\norm{u^m_s}_U}ds,\nonumber
\end{align}
which means $\sup\limits_{0\leq t\leq T; m,n\in \mathbb{N}}\beta^{m,n}(t)< \infty$. So, from $(\ref{4})$ and (\ref{31}), $\{z^n\}$ is a Cauchy sequence in $\mathcal{C}\para{[0,T];H}$, and converges to some $z\in \mathcal{C}\para{[0,T];H}$. Now, we must prove that $z$ is a solution of $(\ref{15})$ with the control function $u$.
Since $\norm{z^n-z}_\infty \rightarrow 0$, $f$ is demicontinuous, and $S(t)$ is contraction, it follows that $S(t-s)f(z^n_s)$ converges to $S(t-s)f(z_s)$ in the weak topology of $H$ for each $0\leq s\leq t\leq T$. Moreover, according to the above a priori estimate and Hypothesis 2(i, iii), we get $\sup\limits_{0\leq s\leq t, m,n\in \mathbb{N}}\norm{S(t-s)f(z^n_s)}<\infty$. Therefore, the Lebesgue dominated convergence theorem implies that, for every $v\in L^2([0,T];H)$,
\begin{align}
\displaystyle\int^t_0\inp{S(t-s)\para{f(z^n_s)-f(z_s)},v(s)}ds\rightarrow0.\label{30}
\end{align}
Since the family $\{S(t-\cdot)\para{f(z^n_{\cdot})-f(z_{\cdot})}: n\in \mathbb{N}\}$ is bounded in the strong topology of $L^2\para{[0,t];H}$, it has an accumulation point in this space with the weak topology, and thanks to (\ref{30}) the accumulation point has to be zero. So, we have for every $t\in [0,T]$
\begin{align}
\norm{S(t-\cdot)\para{f(z^n_\cdot)-f(z_\cdot)}}_{L^2\para{[0,t];H}}\rightarrow0\label{29},
\end{align}
and consequently, $\norm{S(t-\cdot)\para{f(z_\cdot)-f(z^n_\cdot)}}_{L^1\para{[0,t];H}}\rightarrow0$. This means that
$$\displaystyle\int^t_0S(t-s)f(z^n_s)ds\rightarrow\displaystyle\int^t_0S(t-s)f(z_s)ds,$$
for every $t\in [0,T]$. On the other hand
\begin{align}
&\int^t_0S(t-s)\left(g(z^n_s)u^n_s-g(z_s)u_s\right)ds\nonumber\\
&\leq\int^t_0\para{\norm{g(z^n_s)}_L\norm{u^n_s-u_s}_U+D\norm{z^n_s-z_s}\norm{u_s}_U}ds.\label{40}
\end{align}
Besides, from (\ref{10}) and Remark 2.1(ii), we get $\sup\limits_{0\leq s\leq t}\norm{g(z^n_s)}_L<\infty$. Thus, since
$$\norm{u^n-u}_{L^2([0,T];U)}\rightarrow 0\quad\mbox{and}\quad \norm{z^n-z}_\infty \rightarrow 0,$$
we conclude that the right hand side of (\ref{40}) tends to zero. This means that for every $t\in [0,T]$
\begin{align}
\displaystyle\int^t_0S(t-s)g(z^n_s)u^n_sds \rightarrow \displaystyle\int^t_0S(t-s)g(z_s)u_sds.\nonumber
\end{align}
Hence, we obtain
$$z(t)=S(t)x_0+\displaystyle\int^t_0S(t-s)f(z_s)ds+\displaystyle\int^t_0S(t-s)g(z_s)u_sds,$$
for every $t\in [0,T]$ and this completes the proof.
\end{proof}
\end{thm}
The following estimate is an immediate consequence of Remark 2.1 and the a priori estimate $(\ref{10})$.
\begin{cor}\label{c1}
With the same notations as in the proof of Theorem \ref{h1}, we have
\begin{align}
\sup\limits_{n\in \mathbb{N}}\sup\limits_{0\leq t\leq T} \{\norm{f(z^n_t)}  , \norm{g(z^n_t)}  \} <\infty.\label{11}
\end{align}
\end{cor}

\section{\textbf{The Action functional}}
In this section, we deal with Hypothesis 1(i) which concerns the rate function or the action functional. Basically, the rate function is given by ($\ref{22}$), in which $\mathcal{G}^0\para{\int^._0 u(s)ds}$ is the solution of $(\ref{15})$, for any $u\in L^2([0,T];U)$. In the following theorem we prove that the level sets of this rate function are compact which means that we have a good rate function.
\begin{thm}\label{h2}
For every $0<N<\infty$, the set
$$K_N=\set{\mathcal{G}^0\para{\displaystyle\int^._0u_sds}:u\in \mathcal{S}_N}$$
is a compact subset of $\mathcal{C}\para{[0,T];H}$.

\begin{proof} It is sufficient to prove that for every sequence $\{u^n\}$ in $\mathcal{S}_N$, the sequence $\left\{z^n=\mathcal{G}^0\para{\int^._0u^n_sds}\right\}$ has a convergent subsequence in $\mathcal{C}\para{[0,T];H}$. Since $\mathcal{S}_N$ is weakly compact, the sequence $\{u^n\}$ has a convergent subsequence in the weak topology of $\mathcal{S}_N$, still denoted by $\{u^n\}$. We consider $u^0$ as its limit and  show that the sequence $\{z^n\}$ converges to  $z^0=\mathcal{G}^0\para{\int^._0u^0_sds}$ in the space $\mathcal{C}([0,T];H)$.
Let $A_k=A(I-k^{-1}A)^{-1}, k\in \mathbb{N}$, be the Yosida approximation of the operator $A$, and consider the following equation for every $n\geq 0, k>0$:
$$\frac{d}{dt}z^{n,k}_t=A_kz^{n,k}+f(z^{n,k}_t)+g(z^{n,k}_t)u^n_t, \quad z(0)=X_0.$$
From the boundedness of $A_k$, the mild solutions of these equations obtained according to Theorem \ref{h1}, are strong solutions and can be written as
\begin{align}
z^{n,k}_t=\emph{e}^{tA_k}X_0+\displaystyle\int^t_0\emph{e}^{(t-s)A_k}\para{f(z^{n,k}_s)+g(z^{n,k}_s)u^n_s}ds\nonumber.
\end{align}
We show the convergence $\norm{z^n-z^0}_\infty\rightarrow 0$ in two steps. In the first step, we prove that $\norm{z^n-z^{n,k}}_{\infty}\rightarrow 0$, as $k\rightarrow\infty$, for any fixed $n\geq 0$. Then, in the second step, we show that the convergence $\norm{z^{n,k}-z^{0,k}}_\infty\rightarrow 0$ is uniform with respect to $k\in \mathbb{N}$, as $n\rightarrow \infty$. So, according to the first step for $n=0$, we have for every $\eps>0$
$$ \exists k_1\in\mathbb{N}\;\; s.t.\;\;\forall k\geq k_1\quad\norm{z^0-z^{0,k}}_{\infty}<\frac{\eps}{3},$$
and using the uniform convergence from the second step
$$\exists n_1\in\mathbb{N}\;\; s.t.\;\;\forall n\geq n_1,\; k\in\mathbb{N}\quad \norm{z^{n,k}-z^{0,k}}_\infty<\frac{\eps}{3},$$
and finally, according to the first step for $n=n_1$
$$ \exists k_2\in\mathbb{N}\;\; s.t.\;\;\forall k\geq k_2\quad \norm{z^{n_1}-z^{n_1,k}}_{\infty}<\frac{\eps}{3}.$$
Therefore, if we take $k=\max \{k_1, k_2\}$, then $$\norm{z^{n_1}-z^0}_{\infty}<\norm{z^{n_1}-z^{n_1,k}}_{\infty}+\norm{z^{n,k}-z^{0,k}}_\infty+\norm{z^0-z^{0,k}}_{\infty}<\eps.$$\

\emph{Step 1.} By definitions, we have
\begin{align}
z^n_t-z^{n,k}_t=&\displaystyle\int^t_0\para{S(t-s)-e^{(t-s)A_k}}\para{f(z^n_s)+g(z^n_s)u^n_s}ds\nonumber\\
&+\displaystyle\int^t_0 e^{(t-s)A_k}\left(f(z^n_s)-f(z^{n,k}_s)+\para{g(z^n_s)-g(z^{n,k}_s)}u^n_s\right)ds\nonumber\\
=:&I^{n,k}_t+J^{n,k}_t.\label{1}
\end{align}
First, we estimate the norm $\norm{I^{n,k}}_\infty$. Defining, for every $0\leq s\leq T$ and $k\in \mathbb{N}$,
$$\gamma^k(s)=\sup\limits_{s\leq t \leq T}\norm{(S(t-s)-e^{(t-s)A_k})\left(f(z^n_s)+g(z^n_s)u^n_s\right)},$$
we have for every $t\in [0,T]$
\begin{align}
\norm{I^{n,k}_t}\leq \displaystyle\int^t_0\norm{(S(t-s)-e^{(t-s)A_k})\para{f(z^n_s)+g(z^n_s)u^n_s}}ds\leq \displaystyle\int^T_0\gamma^k(s)ds.\nonumber
\end{align}
At the same time, since $\sup\limits_{0\leq t\leq T,k\in \mathbb{N}}\norm{S(t)-e^{tA_k}}\leq 1+1<\infty$ and $f(z^n_\cdot)+g(z^n_\cdot)u^n_\cdot\in L^1([0,T];H)$, we get $\sup\limits_{k\in \mathbb{N}}\gamma^k\in L^1([0,T];\mathbb{R})$. Moreover, from Yosida approximation properties, $\lim_{k\rightarrow\infty}\gamma^k(s)=0$, for every $0\leq s\leq T$. Thus, the Lebesgue dominated convergence theorem implies that
$$\lim_{k\rightarrow\infty} \displaystyle\int^T_0\gamma^k(s)ds=\displaystyle\int^T_0\lim_{k\rightarrow \infty}\gamma^k(s)ds=0.$$
This implies that
\begin{align}
\norm{I^{n,k}}_\infty\rightarrow0,\quad as\; k\rightarrow \infty.\label{21}
\end{align}
Now, we estimate the norm $\norm{J^{n,k}_t}$. We use the energy-type inequality, Hypothesis 2(iv), and Remark
2.1(ii) to write
\begin{align}
\norm{J^{n,k}_t}^2=&\norm{-I^{n,k}_t+z^n_t-z^{n,k}_t}^2\nonumber\\
\leq&\displaystyle\int^t_0 \inp{-I^{n,k}_t+z^n_s-z^{n,k}_s,f(z^n_s)-f(z^{n,k}_s)+\para{g(z^n_s)-g(z^{n,k}_s)}u^n_s}ds\nonumber\\
\leq&\para{\int^t_0\norm{I^{n,k}_s}^2 ds}^{\frac{1}{2}}\para{\int^t_0\norm{f(z^n_s)-f(z^{n,k}_s)}^2 ds}^{\frac{1}{2}}\nonumber\\
&+\para{\int^t_0\norm{I^{n,k}_s}^2 ds}^{\frac{1}{2}}\para{\int^t_0\norm{\para{g(z^n_s)-g(z^{n,k}_s)}u^n_s}^2 ds}^{\frac{1}{2}}\nonumber\\
&+\int^t_0\para{M+D\norm{u^n_s}_U}\norm{z^n_s-z^{n,k}_s}^2ds.\label{17}
\end{align}
Thanks to an estimate analogous to (\ref{11}) for $z^{n,k}$ and the
convergence $(\ref{21})$, the first and second
terms tend to zero, as $k\rightarrow\infty$. Thus, for any
$\delta >0$, one can find large enough $k$ such that
\begin{align}
\norm{J^{n,k}_t}^2\leq\displaystyle\int^t_0\para{M+D\norm{u^n_s}_U}\norm{z^n_s-z^{n,k}_s}^2ds+\delta.\nonumber
\end{align}
So, thanks to $(\ref{21})$ and $(\ref{1})$, we have for a large enough $k$
\begin{align}
\norm{z^n_t-z^{n,k}_t}^2\leq\displaystyle\int^t_0\para{M+D\norm{u^n_s}_U}\norm{z^n_s-z^{n,k}_s}^2ds+\delta.\nonumber
\end{align}
Then, using Gronwall's lemma and a suitable constant $C$, we get
\begin{align}
\norm{z^n_t-z^{n,k}_t}^2\leq C\delta\exp\para{\displaystyle\int^t_0\norm{u^n_s}_Uds}.\label{23}
\end{align}
Since $\delta>0$ is arbitrary, we conclude  for every $n\geq 0$
that
\begin{align}
\lim_{k\rightarrow\infty}\norm{z^n-z^{n,k}}_\infty=0.\label{27}
\end{align}

\emph{Step 2.} By energy-type inequality and due to the fact that $A_k$ is negative definite, we get
\begin{align}
\norm{z^{n,k}_t-z^{0,k}_t}^2 =& \;2\displaystyle\int^t_0\inp{z^{n,k}_s-z^{0,k}_s,A_kz^{n,k}_s-A_kz^{0,k}_s}ds\nonumber\\
&+2\displaystyle\int^t_0\inp{z^{n,k}_s-z^{0,k}_s,f(z^{n,k}_s)-f(z^{0,k}_s)}ds\nonumber\\
&+ 2\displaystyle\int^t_0\inp{z^{n,k}_s-z^{0,k}_s,\para{g(z^{n,k}_s)-g(z^{0,k}_s)}u^n_s}ds\nonumber\\
&+2\displaystyle\int^t_0\inp{z^{n,k}_s-z^{0,k}_s,g(z^{0,k}_s)(u^n_s-u^0_s)}ds\nonumber\\ \leq&2\displaystyle\int^t_0\para{M+D\norm{u^n_s}_U}\norm{z^{n,k}_s-z^{0,k}_s}^2 ds\nonumber\\
&+2\displaystyle\int^t_0\inp{z^{n,k}_s-z^{0,k}_s,g(z^{0,k}_s)(u^n_s-u^0_s)}ds.\label{6}
\end{align}
Estimation of the last term is the challenging problem of this step. For this purpose, we first define
\begin{align}
h^{n,k}(t):=\int^t_0g(z^{0,k}_s)(u^n_s-u^0_s)ds.\nonumber
\end{align}
Since $\{u^n\}$ converges weakly to $u^0$,
$\{g(z^{0,k})(u^n-u^0)\}$ converges weakly to zero in the space
$L^2([0,T];H)$, as $n\rightarrow\infty$. Then, from Lemma
$\ref{lem2}$, it follows that
$\lim_{n\rightarrow\infty}\norm{h^{n,k}}_{\infty}=0$, for any
$k$. Besides, thanks to (\ref{27}), we have
\begin{align}
\lim_{k\rightarrow\infty}\sup\limits_{0\leq s\leq T}\norm{g(z^0_s)-g(z^{0,k}_s)}_L=0,\nonumber
\end{align}
and consequently
\begin{align}
h^{n,k}\rightarrow h^n:=\int^t_0g(z^0_s)(u^n_s-u^0_s)ds,\nonumber
\end{align}
in the space $\mathcal{C}\para{[0,T];H}$, as $k\rightarrow
\infty$. Moreover, the convergence is uniform with respect to
$n$, as $u^n\in \mathcal{S}_N$. This implies that the convergence
$\lim_{n\rightarrow \infty}\norm{h^{n,k}}_{\infty}=0$ is uniform
with respect to $k\in \mathbb{N}$, or equivalently
\begin{align}
\lim_{n\rightarrow\infty}\sup\limits_{k\in \mathbb{N}}\norm{h^{n,k}}_{\infty}=0.\label{34}
\end{align}
Now, using the integration by parts formula for the strong solution, we get
\begin{align}
&\displaystyle\int^t_0\inp{z^{n,k}_s-z^{0,k}_s,g(z^{0,k}_s)(u^{n}_s-u^0_s)}ds\nonumber\\
&=\inp{z^{n,k}_t-z^{0,k}_t,h^{n,k}_t}-\displaystyle\int^t_0\inp{(z^{n,k}_s-z^{0,k}_s)^\prime,h^{n,k}_s}ds\nonumber\\
&\leq\norm{z^{n,k}_t-z^{0,k}_t} \norm{h^{n,k}_t} -\displaystyle\int^t_0\inp{(z^{n,k}_s-z^{0,k}_s)^\prime,h^{n,k}_s}ds.\label{7}
\end{align}
On the other hand,
\begin{align}
&\int^t_0\inp{(z^{n,k}_s-z^{0,k}_s)^\prime,h^{n,k}_s}ds\nonumber\\
&=\int^t_0\inp{A_k({z^{n,k}_s-z^{0,k}_s})+f(z^{n,k}_s)-f(z^{0,k}_s)+g(z^{n,k}_s)u^{n}_s-g(z^{0,k}_s)u^0_s,h^{n,k}_s}ds\nonumber\\
&\leq \norm{h^{n,k}}_{\infty}\displaystyle\int^t_02k\norm{z^{n,k}_s-z^{0,k}_s}ds\nonumber\\
&\quad+\norm{h^{n,k}}_{\infty}\displaystyle\int^t_0\para{\norm{f(z^{n,k}_s)}+\norm{f(z^{0,k}_s)}+\norm{g(z^{n,k}_s)u^{n}_s}+
\norm{g(z^{0,k}_s)u^0_s}ds}.\label{8}
\end{align}
Therefore, from (\ref{34}), we conclude that $(\ref{8})$ and consequently
$(\ref{7})$ tends to zero, uniformly with respect to $k$, as
$n\rightarrow \infty$. Now, if we apply Gronwall's lemma to
$(\ref{6})$, we obtain
\begin{align}
\lim_{n\rightarrow\infty}\sup\limits_{k\in \mathbb{N}}\norm{z^{n,k}-z^{0,k}}_\infty=0.\label{28}
\end{align}
The proof is complete.
\end{proof}
\end{thm}

\section{\textbf{Large deviation principle}}

In this section we deal with the main part of Hypothesis 1. For
every cylindrical Brownian motion $W$ and $\eps>0$, suppose
that $\mathcal{G}^\eps(\eps W)=X^\eps$, in which $X^\eps$ is the mild solution of
$(\ref{13})$. Although in our case the maps $\mathcal{G}^\eps$
are the same for different $\eps>0$, we use the indices
$\{\eps>0\}$ for the sake of simplicity in using Theorem
$\ref{h5}$. The measurability of the maps
$\mathcal{G}^0$ and $\mathcal{G}^\eps$, $\eps>0$, which has
been assumed implicitly in Hypothesis 1, is
a direct consequence of Yamada-Watanabe theorem for mild
solutions. This theorem has been proved by M. Ondrej\'{a}t in \cite{kn:On}.

By Girsanov's theorem and uniqueness of the mild
solution, it follows that $z^\eps_t=\mathcal{G}^\eps\para{\eps
W+\int^._0u^\eps(s)ds}$ is the solution of
\begin{align}
z^\eps_t=S(t)X_0+\displaystyle\int^t_0S(t-s)\para{f(z^\eps_s)+g(z^\eps_s)u^\eps_s}ds+\eps\displaystyle\int^t_0S(t-s)g(z^\eps_s)dW_s.\label{32}
\end{align}
However, for a general $u^{\eps}\in L^2([0,T];U)$, the
nonlinearity of this equation is not semimonotone, uniformly with
respect to $t$. Thus, we can not claim, from Proposition \ref{p2}
or even Theorem 4 in \cite{kn:Za}, that the  solution for this
equation exists. But, similar to Theorem \ref{h1}, if we
approximate a general $u^{\eps}\in L^2([0,T];U)$ by a sequence
$u^{\eps,n}\in L^{\infty}([0,T];U), n\in \mathbb{N}$, we can
overcome this problem easily and find an a priori estimate for the
solution.
\begin{prop}\label{p3}
Let $\mathbb{E} \norm{X_0}^2<\infty$. If the semigroup $S(t)$ and
the functions $f, g$ satisfy Hypothesis 2 and $u^{\eps}\in
\mathcal{P}_N(U)$, for every $\eps>0$, then (\ref{32}) has a
continuous adapted strong solution $z^{\eps}$ such that
\begin{align}
\sup\limits_{0<\eps\leq 1}\mathbb{E}\para{\sup\limits_{0\leq s\leq t}\norm{z^{\eps}_s}^2}<\infty,\quad\forall t\in [0,T].\label{35}
\end{align}
\end{prop}
To verify Hypothesis 1(ii), we must estimate $\norm{z^\eps-z^0}_\infty$ in distribution.
\begin{thm}\label{h3}
Let $\{u^\eps:\eps>0\} \subseteq \mathcal{P}_N(U)$ for some
$N<\infty$. Assume that $u^\eps$ converges to $u^0$ in
distribution and in the weak topology (as
$\mathcal{S}_N(U)$-valued random variables), as $\eps\rightarrow
0$. Then
$$z^{\eps}=\mathcal{G}^\eps\para{\eps W+\displaystyle\int^._0u^\eps_sds}\rightarrow z^0=\mathcal{G}^0\para{\displaystyle\int^._0u^0_sds}$$
in distribution in the space $\mathcal{C}([0,T];H)$.
\begin{proof}
The proof is similar to that of Theorem $\ref{h2}$. The
differences are in some stochastic integrals for which we will
use Lemma $\ref{lem1}$ and Burkholder-type inequality to simplify
them. Let, as before, $A_k$, $k\in \mathbb{N}$, be the Yosida
approximation of the linear operator $A$. Since $A_k$ is bounded,
the equation
$$dz^{\eps,k}_t=\para{A_kz^{\eps,k}_t+f(z^{\eps,k}_t)+g(z^{\eps,k}_t)u^\eps_t}dt+\eps g(z^{\eps,k}_t)dW_t,\quad z(0)=X_0,$$
has a strong solution which can be represented in the form
$$z^{\eps,k}_t=e^{(t-s)A_k}X_0+\int^t_0e^{(t-s)A_k}\para{f(z^{\eps,k}_s)+g(z^{\eps,k}_s)u^\eps_s}ds+\eps\int^t_0e^{(t-s)A_k}g(z^{\eps,k}_s)dW_s.$$
By the definitions, we have
\begin{align}
&z^\eps_t-z^{\eps,k}_t\nonumber\\
=&\displaystyle\int^t_0\para{S(t-s)-\emph{e}^{(t-s)A_k}}\left(f(z^\eps_s)+g(z^\eps_s)u^\eps_s\right)ds\nonumber\\
&+\eps \displaystyle\int^t_0\para{S(t-s)-\emph{e}^{(t-s)A_k}}g(z^\eps_s)dW_s\nonumber\\
&+\displaystyle\int^t_0\emph{e}^{(t-s)A_k}\para{f(z^\eps_s)-f(z^{\eps,k}_s)+g(z^\eps_s)u^\eps_s-g(z^{\eps,k}_s)u^\eps_s}ds\nonumber\\
&+\eps \displaystyle\int^t_0\emph{e}^{(t-s)A_k}\para{g(z^\eps_s)-g(z^{\eps,k}_s)}dW_s=:J^{\eps,k}_1(t)+J^{\eps,k}_2(t)+J^{\eps,k}_3(t)+J^{\eps,k}_4(t).\nonumber
\end{align}
Applying the Burkholder-type inequality, the second and fourth
integrals are estimated respectively by
$$\mathbb{E}\; \norm{J^{\eps,k}_2}^2_{\infty}\leq \eps \mathbb{E}\int^T_0 \norm{g(z^\eps_t)}^2_Ldt\leq C\eps$$
and
$$\mathbb{E} \norm{J^{\eps,k}_4}^2_{\infty}\leq \eps D \mathbb{E}\int^T_0 \norm{z^{\eps,k}_t-z^\eps_t}^2dt\leq C\eps.$$
Similar to the argument used in the previous section to estimate
$\norm{I^{n,k}}_{\infty}$ and $\norm{J^{n,k}}_{\infty}$,  we
conclude that
$\norm{J^{\eps,k}_1}_{\infty},\norm{J^{\eps,k}_3}_{\infty}\rightarrow0$
in distribution, as $k\rightarrow\infty$. Hence, we have the
following convergence in distribution for any $\eps>0$:
\begin{align}
\lim_{k\rightarrow\infty}\norm{z^{\eps}-z^{\eps,k}}_{\infty}\rightarrow0.\label{36}
\end{align}

We now proceed in another direction and let $\eps\rightarrow 0$.
Using It\^{o}'s formula for the strong solution and since $A_k$ is
negative definite, we have
\begin{align}
&\norm{z^{\eps,k}_t-z^{0,k}_t}^2 \nonumber\\
&\quad\leq 2\int^t_0\inp{z^{\eps,k}_s-z^{0,k}_s,A_k\para{z^{\eps,k}_s-z^{0,k}_s}+f(z^{\eps,k}_s)-f(z^{0,k}_s)}ds\nonumber\\
&\qquad+2\int^t_0\inp{z^{\eps,k}_s-z^{0,k}_s,\para{g(z^{\eps,k}_s)-g(z^{0,k}_s)}u^{\eps}_s+g(z^{0,k}_s)(u^\eps_s-u^0_s)}ds\nonumber\\
&\qquad+2\eps\int^t_0\inp{z^{\eps,k}_s-z^{0,k}_s,g(z^{\eps,k}_s)dW_s}+
\eps\int^t_0\norm{g(z^{\eps,k}_s)}^2_Lds\nonumber\\
&\quad\leq 2\int^t_0\para{M+D\norm{u^{\eps}_s}_U}\norm{z^{\eps,k}_t-z^{0,k}_t}^2ds\nonumber\\
&\qquad+2\int^t_0\inp{z^{\eps,k}_s-z^{0,k}_s,g(z^{0,k}_s)(u^\eps_s-u^0_s)}ds\nonumber\\
&\qquad+2\eps\int^t_0\inp{z^{\eps,k}_s-z^{0,k}_s,g(z^{\eps,k}_s)dW_s}+
\eps\int^t_0\norm{g(z^{\eps,k}_s)}^2_Lds\nonumber\\
&\quad=:I^{\eps,k}_1(t)+I^{\eps,k}_2(t)+I^{\eps,k}_3(t)+I^{\eps,k}_4(t).\label{12}
\end{align}
To estimate $I^{\eps,k}_2$, we use the integration by parts
formula for the strong solution. First we define
$$h^{\eps,k}(t)=\int^t_0 g(z^{0,k}_s)(u^\eps_s-u^0_s)ds.$$
From Corollary 3.2 and the fact that $u^{\eps}\rightarrow u^0$ in
distribution and weakly in $L^2([0,T];U)$, it follows that
$\{g(z^{0,k})(u^{\eps}-u^0)\}$ converges to zero in distribution
and weakly in the space $L^2([0,T];H)$, as $\eps\rightarrow0$.
Due to Lemma $\ref{lem1}$, this implies that
$\lim_{\eps\rightarrow0}\norm{h^{\eps,k}}_\infty=0$ in
distribution, for any $k\in\mathbb{N}$. Furthermore, thanks to
(\ref{27}), we conclude that
\begin{align}
h^{\eps,k}\rightarrow h^{\eps}:=\int^t_0 g(z^0_s)(u^\eps_s-u^0_s)ds,\nonumber
\end{align}
in the space $\mathcal{C}([0,T];H)$ and in distribution, as $k\rightarrow\infty$. So,
we have in distribution
\begin{align}
\lim_{\eps\rightarrow0}\sup\limits_{k\in \mathbb{N}}\norm{h^{\eps,k}}_{\infty}=0.\label{37}
\end{align}
Therefore, using the integration by parts for the strong solution, we get
\begin{align}
&\frac{1}{2}I^{\eps,k}_2(t)=\inp{z^{\eps,k}_t-z^{0,k}_t,h^{\eps,k}_t} \nonumber\\
&-\displaystyle\int^t_0\inp{A_k({z^{\eps,k}_s-z^{0,k}_s}),h^\eps_s}ds-\displaystyle\int^t_0\inp{f(z^{\eps,k}_s)-f(z^{0,k}_s),h^\eps_s}ds\nonumber\\
&-\displaystyle\int^t_0\inp{g(z^{\eps,k}_s)u^\eps_s-g(z^{0,k}_s)u^0_s,h^\eps_s}ds-\eps\displaystyle\int^t_0\inp{g(z^{\eps,k}_s)dW_s,h^\eps_s}.\label{19}
\end{align}
Applying  Lemma $\ref{lem1}$ to the last term and the
Cauchy-Schwarz inequality to the other terms, we obtain
\begin{align}
\mathbb{E}\norm{I^{\eps,k}_2}_{\infty} \leq&2\mathbb{E}\left\{\norm{h^\eps}_{\infty}\norm{z^{\eps,k}-z^{0,k}}_{\infty}\right\}+4k\mathbb{E}\left\{\norm{h^{\eps,k}}_{\infty}\int^T_0\norm{z^{\eps,k}_t-z^{0,k}_t}dt\right\}\nonumber\\
&+2\mathbb{E}\left\{\norm{h^{\eps,k}}_{\infty}\int^T_0 \para{\norm{f(z^{\eps,k}_t)}+\norm{f(z^{0,k}_t)}}dt\right\}\nonumber\\
&+2\mathbb{E}\left\{\norm{h^{\eps,k}}_{\infty}\int^T_0 \para{\norm{g(z^{\eps,k}_t)u^\eps_t}+\norm{g(z^{0,k}_t)u^0_t}}dt\right\}\nonumber\\
&+3\eps\mathbb{E}\left\{\norm{h^{\eps,k}}_{\infty}+\int^T_0\norm{g(z^{\eps,k}_t)}^2_Ldt\right\}.\nonumber
\end{align}
So, from (\ref{37}) and the estimate (\ref{35}) for
$z^{\eps,k}$, we get in distribution $$\lim_{\eps\rightarrow 0}\sup\limits_{k\in \mathbb{N}}\norm{I^{\eps,k}_2}_{\infty}=0$$.

To estimate the term $I^{\eps,k}_3$, we use Lemma $\ref{lem1}$. So we have
\begin{align}
\mathbb{E}\norm{I^{\eps,k}_3}_{\infty}\leq 3\eps\mathbb{E}\left\{\norm{z^{\eps,k}-z^{0,k}}^2_{\infty}+\int^T_0\norm{g(z^{\eps,k}_t)}^2_Ldt\right\}.\label{24}
\end{align}
Then, by the a priori estimate (\ref{35}), we obtain $\lim_{\eps\rightarrow 0}\sup\limits_{k\in \mathbb{N}}\norm{I^{\eps,k}_3}_{\infty}=0$ in distribution. Again using (\ref{35}), it is clear that $\lim_{\eps\rightarrow 0}\sup\limits_{k\in \mathbb{N}}\norm{I^{\eps,k}_4}_{\infty}=0$ in distribution. Therefore, applying Gronwall's lemma to $(\ref{12})$, we obtain the following convergence in distribution:
\begin{align}
\lim_{\eps\rightarrow 0}\sup\limits_{k\in \mathbb{N}}\norm{z^{\eps,k}-z^{0,k}}_\infty=0.\nonumber
\end{align}
Finally, by an argument similar to that of Theorem 4.1, we conclude that\\
$\lim_{\eps\rightarrow 0}\norm{z^{\eps}-z^0}_{\infty}=0$ in distribution and the proof is complete.
\end{proof}
\end{thm}

\section{\textbf{Examples}}
\subsection{\textbf{Stochastic heat equation}}
As the first example, we consider the LDP for the semilinear stochastic heat equation with Dirichlet boundary condition and monotone nonlinear term. This example is a direct generalization of the problem considered in \cite{kn:FJ}. Faris and Jona-Lasinio  studied in \cite{kn:FJ} the LDP and the tunneling problem for the heat equation with additive noise and nonlinear term $f(x)=\lambda x^4-\mu x^2, \; \lambda , \mu >0$. The LDP for the heat equation with Lipschitz nonlinearity and multiplicative noise has been investigated  by Peszat in \cite{kn:Pe}.
\begin{defn}
An $\mathbb{R}^m$-valued function $f(x,u)$ of two variables $x\in D\subseteq \mathbb{R}^d, u\in \mathbb{R}^n$ is said to satisfy the Carath\'{e}odory condition, if it is continuous with respect to $u$ for almost all $x\in D$ and measurable with respect to $x$ for all values of $u\in \mathbb{R}^n$.
\end{defn}
Let $D$ be a bounded domain with smooth boundary in $\mathbb{R}^d$. Consider the initial-boundary value problem
\begin{align}
\begin{array}{c}
  \frac{\partial u}{\partial t}=\Delta u+f(x,u(t,x))+\Sigma^{l}_{i=1}g_i(x,u(t,x))\dot{W}_i(t), \quad (t,x)\in [0,\infty)\times D,\\\label{25}
  u(t,x)=0, \quad \forall (t,x)\in [0,\infty)\times \partial D,\\
  u(0,x)=u_0(x) \quad \forall x\in D;
\end{array}
\end{align}
in which, $\Delta$ is the Laplacian operator, $W_i$, $i=1,\ldots,l$, are independent standard
real Brownian motions and $u_0\in L^2(D)$. Moreover, $f$ and $g_i$, $i=1,\ldots,l$, satisfy the following hypothesis.\\\\
\textbf{Hypothesis 3}
\begin{description}
  \item (i) $f,g_i:D\times \mathbb{R}\rightarrow\mathbb{R}$, $i=1,\ldots,l$, satisfy the Carath\'{e}odory condition;
  \item (ii) There exist a function $a\in L^2(D)$ and a constant $C>0$ such that for each $0\leq i\leq l$
$$\abs{g_i(x,y)}+\abs{f(x,y)}\leq a(x)+C\abs{y}, \quad \forall \; (x,y)\in D\times \mathbb{R};$$
  \item (iii) $g_i(x,.)$, for each $0\leq i\leq l,$ is uniformly Lipschitz with Lipschitz constant $C>0$, i.e.
$$\abs{g_i(x,y_2)-g_i(x,y_1)}\leq C\abs{y_2-y_1}, \quad \forall x\in D,\;  y_1, y_2\in \mathbb{R};$$
  \item (iv) $-f(x,.)$ is uniformly semimonotone with parameter M, i.e.
$$\inp{f(x,y_2)-f(x,y_1),y_2-y_1}\leq M(y_2-y_1)^2, \quad \forall x\in D,\;  y_1, y_2\in \mathbb{R}.$$
\end{description}
Let $H=L^2(D)$, with the norm $\norm{\cdot}_{L^2(D)}$, and $A=\Delta u$ with the domain
$$D(A)=\left\{u\in L^2(D)\; :\; u', u''\in L^2(D), \; u(x)=0\; \forall x\in \partial D \right\}.$$
The operator $A:D(A)\subset H\rightarrow H$ generates a
strongly continuous semigroup $S(t), t\geq 0$. Now, define
$\bar{f}, \bar{g_i}:H\rightarrow H$, $0\leq i\leq l$, by
\begin{align}
\begin{array}{c}
  \bar{f}(u)(x)=f(x,u(x)), \quad u\in H, \; x\in D;\\
  \bar{g}_i(u)(x)=g_i(x,u(x)) \quad u\in H, \; x\in D.\nonumber
\end{array}
\end{align}
Since $f$ and $g_i$, $i=1,\ldots,l$, satisfy Hypothesis 3(i, ii), then from \cite[Theorem 2.1]{kn:Kra}, $\bar{f}$ and $\bar{g}_i$, $i=1,\ldots,l$, are continuous and there is $C>0$ such that
$$\norm{\bar{f}(u)}\leq C(1+\norm{u})\;\; and \; \; \norm{\bar{g}_i(u)}\leq C(1+\norm{u}).$$
According to  Hypothesis 3(iii, iv), it is clear that $\bar{f}$ is semimonote and $\bar{g}_i$ satisfies the Lipschitz condition. Define $\bar{g}=(\bar{g}_1,\ldots,\bar{g}_l):H\rightarrow(L^2(D))^l\simeq L(\mathbb{R}^l;L^2(D))$. Note that $L(\mathbb{R}^l;L^2(D))=L_2(\mathbb{R}^l;L^2(D))$, since $\Real^l$ is finite dimension. Therefore, we can rewrite $(\ref{25})$ as the following semilinear stochastic evolution equation:
\begin{align}
du(t)=-Au(t)dt+\bar{f}(u(t))dt+\bar{g}(u(t))dW^t(t), \quad u(0)=u_0,\nonumber
\end{align}
where, $W^t=(W_1,\ldots,W_l)^t$. Now, considering $U=\mathbb{R}^l$ and according to Theorem 5.1, we obtain the LDP for small noise limit of the mild solution to $(\ref{25})$.

\subsection{\textbf{Initial-value problem of semilinear hyperbolic systems}}
Consider the following initial-value problem of a semilinear hyperbolic system:
\begin{align}
\begin{array}{c}
  \frac{\partial u}{\partial t}=\sum^n_{i=1} a_i(x)\frac{\partial u}{\partial x_i}+b(x)u+f(x,u)+g(x,u)\dot{W}, \\
  u(0,x)=u_0, \quad u_0(x)\in L^2(\Real^n;\Real^N), \quad x\in \Real^n,\label{26}
\end{array}
\end{align}
where $\dot{W}$ is an $m$-dimensional Brownian motion, $u=(u_1,\ldots,u_N)^t$ is the unknown
vector, and $a_i(x)$ and $b(x)$ are square matrices of order $N$, for $i=1,\ldots,n$ and
$x\in \Real^n$. In this example we assume the following hypotheses:\\\\
\textbf{Hypothesis 4}
\begin{description}
  \item (i) The matrices $a_i(x)$, $i=1,\ldots,n$ and $x\in \Real^n$, are symmetric;
  \item (ii) The components of $a_i$, $i=1,\ldots,n$, their first order derivatives, and the function $b$ are bounded continuous;
  \item (iii) $f:\Real^n\times \Real^N\rightarrow \Real^N$ and $g:\Real^n\times\Real^N\rightarrow L(\Real^m,\Real^N)$ satisfy the Carath\'{e}odory condition;
  \item (iv) There exist a function $a\in L^2(\Real^n)$ and a constant $C>0$ such that for all $x\in \Real^n$ and $u\in \Real^N$
$$\norm{g(x,u)}_{L(\Real^m,\Real^n)}+\norm{f(x,t)}\leq a(x)+C\norm{u};$$
  \item (v) $-f$ and $g$ are semimonotone and Lipschitz, respectively, in the second variable,
  uniformly with respect to the first variable.
\end{description}
Let $H=L^2(\Real^n;\Real^N)$ and $U=\Real^m$. Define a closed
unbounded operator $A$ on $H$ by
$$Au=\sum^n_{i=1} a_i(x)\frac{\partial u}{\partial x_i}+b(x)u, \quad \forall u\in D(A)\subset H.$$
According to \cite[Theorem 3.51]{kn:Ta}, $A$ is the generator of a
$C_0$-semigroup on $H$. we define $\bar{f}:H\rightarrow H$ by
$\bar{f}(u)(x)=f(x,u(x))$ and $\bar{g}:H\rightarrow
L(U,H)=L_2(U,H)$ by $\bar{g}(u)(x)=g(x,u(x))$ for every $u\in H$
and $x\in\Real^n$. Then $(\ref{26})$ reduces to the semilinear
evolution equation
\begin{align}
du_t=Au_tdt+\bar{f}(u_t)dt+\bar{g}(u_t)dW_t, \quad u(0)=u_0.\nonumber
\end{align}
Similar to the previous example, $\bar{f}$ and $\bar{g}$ are continuous and there exist a constant $C>0$ such that
\begin{align}
\begin{array}{c}
  \norm{\bar{f}(u)}\leq C(1+\norm{u}),\\
  \norm{\bar{g}(u)}_L\leq C(1+\norm{u}).\nonumber
\end{array}
\end{align}
So, from Hypothesis 4, it is clear that $A$, $\bar{f}$ and $\bar{g}$ satisfy Hypothesis 2.
Therefore, the mild solution of SPDE $(\ref{26})$ satisfies LDP.

For more examples of semilinear evolution equations with monotone nonlinearity we refer to \cite{kn:Za}.\\\\
\textbf{\textbf{Acknowledgement}}.  The first named author would like
to thank Professor M. R\"{o}ckner for his valuable suggestions.


\emph{Department of Mathematical Sciences}\\
\emph{Sharif University of Technology}\\
\emph{P.O. Box 11365-9415}\\
\emph{Tehran, Iran}\\
\emph{e-mail: dadashi@mehr.sharif.edu}\\
\emph{e-mail: zangeneh@sharif.edu}
\end{document}